\documentclass[letterpaper, 10 pt, conference]{ieeeconf}

\IEEEoverridecommandlockouts                              % This command is only needed if 
\usepackage[utf8]{inputenc}
\usepackage{cite}
\usepackage{graphicx} 
\usepackage{amsmath}
\usepackage{siunitx}
\usepackage{multirow}
\usepackage{comment} 
\usepackage{authblk}
\usepackage{algorithmicx}
\usepackage[ruled]{algorithm2e}% http://ctan.org/pkg/algorithmicx
\usepackage[table]{xcolor}
\usepackage{lettrine}
\usepackage{url}

\allowdisplaybreaks[4]
\newcommand\mcedit[1]{\textcolor{black}{#1}}
\newcommand\mdsedit[1]{\textcolor{black}{#1}}

\newcommand\mbedit[1]{\textcolor{black}{#1}}
\title{\LARGE \bf
Fast optimal trajectory generation for a tiltwing VTOL aircraft with application to urban air mobility
}

\author{Martin Doff-Sotta$^{\star}$\qquad Mark Cannon$^{\star}$\qquad  Marko Bacic$^{\star, \dagger}$% <-this % stops a space
\thanks{$^{\star}$All authors are with the Control Group, University of Oxford, Parks Road, Oxford OX1 3PJ, United Kingdom. Corresponding author: Martin Doff-Sotta ({\tt\small martin.doff-sotta@eng.ox.ac.uk})}%
\thanks{$^{\dagger}$On part-time secondment from Rolls-Royce plc.}%
}

\begin{document}

\maketitle
\pagestyle{plain}

\begin{abstract}
We solve the minimum-thrust optimal trajectory generation problem for the transition of a tiltwing Vertical Take-Off and Landing (VTOL) aircraft using convex optimisation. The method is based on a change of differential operator that allows us to express the simplified point-mass dynamics along a prescribed path and formulate the original nonlinear problem 
%as a sequence of two convex programs.
\mcedit{in terms of a pair of convex programs}. 
A case study involving the Airbus A$^3$ Vahana VTOL aircraft is considered for forward and backward transitions.  The presented approach provides a fast method to generate a safe optimal transition for a tiltwing VTOL aircraft that can further be leveraged online for control and guidance purposes. \\

\end{abstract}

%%%%%%%%%%%%%%%%%%%%%%%%%%%%%%%%%%%%%%%%%%%%%%%%%%%%%%%%%%%%%%%%%%%%%%%%%%%%%%%%%%%%%%%%%%%%%%%%%%%%%%%%%%%%%%%%%%%%%%%%%%%%%%%%%%%%%%%%%%%
% Martin: replaced all US spelling by UK spelling (e.g. optimization -> optimisation) for consistency with the rest of the thesis 
% Mark: OK!
%%%%%%%%%%%%%%%%%%%%%%%%%%%%%%%%%%%%%%%%%%%%%%%%%%%%%%%%%%%%%%%%%%%%%%%%%%%%%%%%%%%%%%%%%%%%%%%%%%%%%%%%%%%%%%%%%%%%%%%%%%%%%%%%%%%%%%%%%%%
\section{Introduction}
\lettrine{T}{he} increase of congestion and pollution levels in large metropolitan areas has recently incentivized the development of sustainable and green mobility solutions. %As an example, in 2014 in the USA, road traffic has generated more than 11 millions litres of excess fuel and nearly 7 billion person-hours were wasted - the equivalent of about 10,000 people average lifespans \cite{schrank20152015}.
Urban Air Mobility (UAM) has the potential to address these problems by allowing air transport of goods and people \cite{easa}, thus reducing the pressure on urban traffic. \mcedit{It has been estimated that 160,000 air taxis could be in circulation worldwide by 2050, representing a USD 90 billion market \cite{roland-berger}.}
%By 2050, an estimated 160,000 air taxis could be in circulation worldwide, representing a USD 90 billion market \cite{roland-berger}. 

Tiltwing VTOL aircraft, with their capability to take-off and land in restricted spaces and their extended endurance, have recently received a lot of attention in a UAM context. Despite being \mcedit{investigated} %around
since the \mcedit{1950's} \cite{kuhn1959semiempirical, hargraves1961analytical, mccormick1967aerodynamics}, issues related to control and stability, as well as the mechanical complexity associated with the use of conventional internal combustion engines have \mcedit{prevented widespread adoption of this technology}. However, recent advances in battery and electric motor technology, combined with the development of modern control systems architectures, has spurred a new interest for these vehicles \cite{UMich19}. \mcedit{Some recent} prototypes based on tiltwing configurations \mcedit{are}
%were proposed recently: 
the Airbus A$^3$ Vahana, the NASA GL-10 \cite{rothhaar2014nasa}, or Rolls-Royce's eVTOL aircraft \cite{higgins2020aeroacoustic}.  

The transition manoeuvre for a VTOL aircraft is the most critical phase of flight. Although heuristic strategies \mcedit{have been}
%were 
proposed based on smooth scheduling functions of the forward velocity or tilt angle \cite{me}, trajectory optimisation for the transition of VTOL aircraft is still an open problem as it involves determining the optimal combination of thrust and tiltwing angle to minimise an objective while meeting \mcedit{constraints on system states and control inputs},
%state and input constraints, 
which requires solving large-scale nonlinear optimisation problems. 

The problem \mcedit{of determining minimum energy speed profiles}
%was addressed in \cite{Iowa19}
%, where a multiphase optimization problem with fixed arrival time is solved 
for the forward transition manoeuvre of the Airbus A$^3$ Vahana \mcedit{was addressed in \cite{Iowa19}}, considering various phases of flight (cruise, transition, descent). A drawback of the approach is that the transition was assumed to occur instantly, which is \mcedit{unrealistic}
%not reasonable
for such \mcedit{a} vehicle. 
%
%The vehicle is modelled in all three phases of flight considered: cruise, transition, descent.  The energy consumption in all flight phases is  established and a flight profile is searched that minimizes the total energy while satisfying the dynamic and kinematic constraints. A major drawback of the approach is that the tiltwing angle was not included in the optimization and the transition was assumed to occur instantly, which is not reasonable for such vehicle. 
%
In \cite{UMich19}, the trajectory generation problem for take-off is formulated as a constrained optimisation problem \mcedit{and} solved using NASA's OpenMDAO framework and the SNOPT gradient-based optimiser. A case study based on the Airbus A$^3$ Vahana tiltwing aircraft is considered, including aerodynamic models of the wing beyond stall, and the effect of propeller-wing interaction.  
%
%The drive power is minimised considering the wing angle, electric power and flight time as optimization variables with constraints on final altitude and terminal cruise speed. Additional constraints on angle of attack, acceleration and horizontal distance are also considered.  A 5\% energy saving is achieved using the proposed optimization when compared to a baseline situation where the aircraft transitions at a constant optimal climb rate. It is shown that the optimal trajectory involves operating the wing in near-stall conditions, although the benefit of doing so is practically negligible. 
%
%The 
Forward and backward optimal transition manoeuvres at constant altitude are computed in \cite{leo} for a tiltwing aircraft, considering leading-edge fluid injection active flow control and the use of a high-drag device. %The effect of propeller vortex interaction with the wing is modelled and included in the model of the vehicle. The optimization problem is transcribed into a NLP using GPOPS-II, a general-purpose nonlinear optimization software implementing a pseudospectral collocation method at the Legendre-Gauss-Radau points. The NLP is then solved using the  IPOPT solver implementing a primal-dual interior point method.   
All \mcedit{of these approaches consider optimising trajectories}
%All the approaches so far considered solving the nonlinear problem directly 
using general purpose NLP solvers without exploiting potentially 
\mcedit{useful}
%interesting 
structures and simplifications. 
\mcedit{This approach is}
%This can be 
computationally expensive \mcedit{and is not generally suitable for real-time implementation}.
%can be computationally expensive and might be intractable online.

We propose in this paper a \mcedit{convexification of the optimal trajectory generation} 
%convex formulation of the 
problem based on a change of variables \mcedit{inspired by \cite{bobrow}}. 
This translates the nonlinear equations of motion in \mcedit{the} time domain into \mcedit{a set of linear differential equations in the space}
%a linear formulation in space 
domain along a prescribed path. 
\mcedit{We thus derive}
%This permits the development of 
an efficient method based on \mcedit{convex programming}. \mdsedit{Convex programming was successfully employed to solve optimisation problems related to energy management of hybrid-electric aircraft \cite{me2, me3} and spacecraft \cite{me4}.}
%the general purpose convex programming software package CVX \cite{cvx} with solver Mosek \cite{mosek}. 
%
%ADMM was successfully employed to solve convex optimization problems related to energy management of hybrid-electric aircraft \cite{me2, me3}, spacecraft rendezvous problems \cite{le2019fast} or model predictive control of Unmanned Aerial Vehicles \cite{cheng2020semi}. The resulting optimal solution could further be leveraged by other parts of the control architecture, e.g. a model predictive controller.
%
The paper is organized as follows. Section \ref{sec:modeling} introduces \mcedit{the} tiltwing VTOL aircraft dynamic and aerodynamic models. \mcedit{Section \ref{sub:convex_formulation} formulates the continuous time trajectory optimisation problem in terms of a pair of}
%In Section \ref{sub:convex_formulation}, the dynamics are expressed in a form suitable for convex programming and the trajectory optimization problem is formulated as two 
convex problems. \mcedit{These} 
%The obtained optimization problems 
are discretised in Section \ref{sec:discrete} and Section \ref{sec:results} discusses simulation results obtained for a case study based on the Airbus A$^3$ Vahana. 
\mcedit{Section~\ref{sec:conclusion} presents conclusions.}

\section{Modeling}
\label{sec:modeling}
 A simplified longitudinal point-mass model of a tiltwing VTOL aircraft equipped with propellers is developed in this section. In order to account for the effect of the propeller wake on the wing, the flow velocity downstream is augmented by the induced velocity of the propeller. %This allows us to define two new variables: the effective velocity and effective angle of attack. These 
%account for modifications in the flow due to the propeller slipstream, thus impacting lift and drag and ultimately \mcedit{the aircraft dynamics}.
%the dynamics. 
The second-order dynamics of the tilting wing is also presented.
 
\begin{figure}[h]
    \centering
    \includegraphics[width=0.45\textwidth, trim={1.5cm 3.5cm 1.5cm 2cm}, clip]{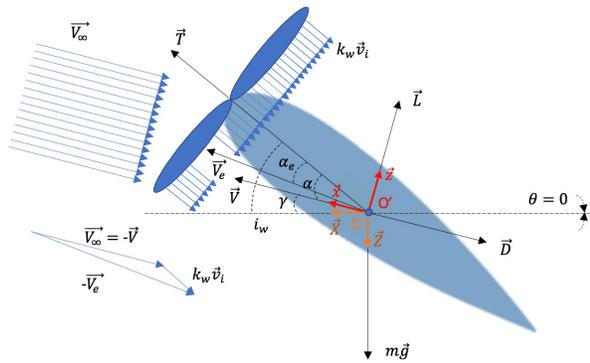} %{<left> <lower> <right> <upper>}
    \caption{\mcedit{Force and velocity definitions for a VTOL aircraft}}
    \label{fig:diagram}
\end{figure}

\mcedit{We}
%Let us 
consider a planar point-mass model of a VTOL aircraft as shown in Figure \ref{fig:diagram}, whose position with respect to an inertial frame $O_{XZ}$ is given by $(x,z)$ and velocity given by %Introducing the body attached frame $O'_{xz}$ whose $x$-axis is aligned with the velocity vector and $z$-axis aligned with the lift vector, we obtain 
%the following parametrization in polar coordinates
%
\[
\dot{x} = V \cos{\gamma},
\qquad
\dot{z} = -V \sin{\gamma},
\]
% \begin{align*}
% \dot{x} &= V \cos{\gamma},
% \\
% \dot{z} &= -V \sin{\gamma},
% \end{align*}
%
where $V$ is the aircraft velocity magnitude and $\gamma$ the \mdsedit{flight path angle, defined as the angle of the velocity vector from horizontal}.
%flight path angle. 
From \mcedit{Figure~\ref{fig:diagram}}
%the diagram
the point-mass equations of motion (EOM) in polar coordinates are
%given below
%
\begin{alignat}{2}
m \dot{V} &= T \cos{\alpha} - D -mg \sin{\gamma}, 
&\quad 
V(t_0) &= V_{0},
\label{eq:eom1}
\\
m V \dot{\gamma} &= T \sin{\alpha} + L -mg \cos{\gamma}, 
&\quad 
\gamma(t_0) &= \gamma_{0},
\label{eq:eom2}
\end{alignat}
where $m$ is the mass of the \mcedit{aircraft},
%vehicle, 
$g$ the acceleration due to gravity,  $T$ the thrust magnitude, $\alpha$ the angle of attack.

The dynamics of the wing are given by 
\begin{equation}
J_w \ddot{i}_w = M, \quad i_w(t_0) = i_{0}, \quad \dot{i}_w(t_0) = \Omega_{0},
\label{eq:tiltwing_dyna}
\end{equation}
where $J_w$ is the rotational inertia of the wing (\mcedit{about}
%along 
the $y$-axis), $M$ is the total torque delivered by the tilting actuators and $i_w$ is the tiltwing angle such that 
\[
i_w + \theta = \alpha + \gamma .
%\label{eq:angle}
\]
\mcedit{Here}
%where 
$\theta$ is the pitch angle, defined as the angle 
\mcedit{of the fuselage axis from horizontal}. 
%between the fuselage axis and the horizon
\mcedit{For passenger comfort}, $\theta$ is regulated via the elevator to track a constant reference $\theta^* = 0$.
%for improved passengers comfort.

From momentum theory, the propeller generates an induced speed $v_i$ such that  
\[
T = \rho A n (V \cos{\alpha} + v_i )( k_w v_i),
\]
%
\begin{comment}
\[
v_i = -\frac{V \cos{\alpha}}{2} + \sqrt{\frac{V^2 \cos^2{\alpha}}{4} + \frac{T}{2\rho A n}},
\]
\end{comment}
%
where $\rho$ is the air density, $A$ the rotor disk area, $n$ the number of \mcedit{propellers, and}
%motors, 
$k_w \approx 2$. The effective velocity $V_e$ and effective angle of attack $\alpha_e$ seen by the wing are given by  
\begin{align*}
V_e \cos{\alpha_e} &= V \cos{\alpha} + k_w v_i ,
%\label{eq:eff2}
\\
V_e \sin{\alpha_e} &= V \sin{\alpha},
%\label{eq:eff2}
\end{align*}
%Combining \mcedit{these} 
%equations yields
\[
V_e^2 = V^2 + \frac{2T}{\rho A n}. 
%\label{eq:eff1}
\]

The total lift and drag are modeled as the weighted sum of the blown and unblown counterparts
%($\alpha, \alpha_e$ small):
%
\begin{align*}
D &= (1-\mu) \tfrac{1}{2}  \rho S(a_2 \alpha^2 + a_1 \alpha  + a_0) V^2  \\
&\quad + \mu \tfrac{1}{2}  \rho S(a_2 \alpha_e^2 +a_1 \alpha_e  + a_0) V_e^2 
\\
&\approx (1-\mu) \tfrac{1}{2}  \rho S(a_1 \alpha  + a_0) V^2  + \mu \tfrac{1}{2}  \rho S(a_1 \alpha_e  + a_0) V_e^2,
%\label{eq:drag}
\\
L &= (1-\mu) \tfrac{1}{2}  \rho S(b_1 \alpha  + b_0) V^2  + \mu \tfrac{1}{2}  \rho S(b_1 \alpha_e  + b_0) V_e^2.
%\label{eq:lift}
\end{align*}
\mcedit{Here}
%where 
$S$ is the wing area, 
$\mu = 2 R n / l$ is the \mdsedit{blown-unblown ratio, $R$ is the radius of the rotor disk, $l$ is the wingspan},
\mcedit{$a_0, a_1, a_2$ and $b_0, b_1$ are constant parameters and 
$\alpha, \alpha_e\ll 1$ is assumed so that terms involving $\alpha^2,\alpha_e^2$ are negligible}. 
The following input and state constraints also apply
\begin{gather*}
 0 \leq T \leq \overline{T}, \quad 0 \leq V \leq \overline{V},
\nonumber
\\
V(t_0) = V_0
\text{ and }
V(t_f) = V_f,
\nonumber
\\
\underline{a} \leq \dot{V} \leq \overline{a}, \quad \underline{M} \leq M \leq \overline{M}.
\nonumber
\end{gather*}

\mcedit{This paper considers how to optimise}
%A current research problem is to generate an optimal trajectory for 
the transition between powered lift and cruise flight modes. 
%of a VTOL aircraft 
%\mcedit{subject to these constraints}.
%whose dynamics is described by the equations above. 
Minimising thrust is a natural choice as a proxy for minimising energy consumption, 
\mcedit{suggesting the following objective function}
\begin{equation}
    \label{eq:obj}
    J = \int_{t_0}^{t_f} \left(T/\overline{T}\right)^2 \mathrm{d}t. 
\end{equation}
The optimisation problem consists of minimising \eqref{eq:obj} while satisfying dynamical constraints, \mcedit{input and state constraints}.

\section{Convex formulation}
\label{sub:convex_formulation}
The \mcedit{aircraft dynamics} 
%equations
derived in Section \ref{sec:modeling} contain 
%some 
nonlinearities \mcedit{that make the}
% which implies that the exact 
trajectory optimisation problem \mcedit{nonconvex}.
%for the presented tiltwing VTOL aircraft is nonlinear.
However, 
%we present here 
a change of \mcedit{variables (motivated by~\cite{bobrow})} 
%variable \cite{bobrow} that can simplify 
considerably \mcedit{simplifies} the structure of \mcedit{this problem}
%the equations above 
and allows us to formulate the problem \mcedit{in terms of}
%as a series of 
two convex programs. If we assume that a path $(x(s), z(s))$ parameterised by the curvilinear abscissa $s$ is already known \textit{a priori} (which \mcedit{is usually}
%would usually be 
the case in a UAM context where flight corridors are prescribed) we can introduce \mcedit{via}
%by 
the chain rule the following change of differential operator,
\begin{equation}
\label{eq:CV}
\frac{\mathrm{d}}{\mathrm{d}t} = V \frac{\mathrm{d}}{\mathrm{d}s} .
\end{equation}
\mcedit{Denoting $\tfrac{\mathrm{d}\,\cdot}{\mathrm{d}s}=\cdot'$, the EOM in \eqref{eq:eom1} and \eqref{eq:eom2} are equivalent to}
\begin{align*}
m V V' &= T \cos{\alpha} - (1-\mu)\frac{1}{2}\rho S (a_1 \alpha + a_0 ) V^2 
\\
&\quad -\mu\frac{1}{2}\rho S (a_1 \alpha_e + a_0) \Bigl(V^2 + \frac{2T}{\rho A n}\Bigr) -mg \sin{\gamma^*},
%\label{eq:eom5}
\\
m V^2 \gamma^{*\prime} &= T \sin{\alpha} + (1-\mu)\frac{1}{2}\rho S (b_1 \alpha + b_0) V^2 
\\
&\quad +\mu\frac{1}{2}\rho S (b_1 \alpha_e + b_0 ) \Bigl(V^2 + \frac{2T}{\rho A n}\Bigr)  -mg \cos{\gamma^*},
%\label{eq:eom6}
\end{align*}
where $\gamma^*$ is known \textit{a priori} from the path.
\mcedit{Specifically},
since $\mathrm{d}x= V \cos{\gamma} \, \mathrm{d}t$ \mcedit{and} $\mathrm{d}z=-V \sin{\gamma} \, \mathrm{d}t$, the flight path angle can be expressed in terms of the path variables \mcedit{as}
\[
\tan{\gamma^*} = -\frac{\mathrm{d}z}{\mathrm{d}x} .
\]

\mcedit{Defining new model states}
%Now introducing the new states 
$E=V^2$ and $a=VV'$, \mcedit{so that}
%and noting that 
%
\begin{equation}
\label{eq:diff}
E' = 2 a, 
\end{equation}
the EOM reduce to 
\begin{align}
\label{eq:eom7}
& m a = T \cos{\alpha} - (1-\mu)\frac{1}{2}\rho S\left(a_1 \alpha + a_0\right) E 
\nonumber\\
&\quad- \mu\frac{1}{2}\rho S\left(a_1 \alpha_e + a_0\right)\left(E + \frac{2T}{\rho A n}\right) -mg \sin{\gamma^*},
\\
\label{eq:eom8}
& m E \gamma^{*\prime} = T \sin{\alpha} + (1-\mu)\frac{1}{2}\rho S\left(b_1 \alpha + b_0\right)E 
\nonumber
\\
&\quad + \mu\frac{1}{2}\rho S\left(b_1 \alpha_e + b_0\right)\left(E + \frac{2T}{\rho A n}\right)  -mg \cos{\gamma^*}. 
\end{align}
\mcedit{Let} $\lambda = a_1/b_1$, \mcedit{then}  
the linear combination
%of these equations above 
\eqref{eq:eom7} $+$ $\lambda$\eqref{eq:eom8} 
%$=$ \eqref{eq:nocvx} 
yields
\begin{multline}
\!\!\!\!\! m a +  \underbrace{(\lambda m  \gamma^{*\prime} + \frac{1}{2}\rho S (a_0 -\lambda b_0))}_{c(\gamma^{*\prime})}E  + \underbrace{mg( \sin{\gamma^*} + \lambda \cos{\gamma^*})}_{d(\gamma^*)} 
\\
= \underbrace{T \cos{\alpha} + \lambda T \sin{\alpha} - \mu S^\star (a_0 - \lambda b_0 ) T}_{\tau} ,
\label{eq:nocvx}
\end{multline}
where 
$S^\star=\frac{S}{AN}$ 
\mcedit{and $\tau$  is a virtual input defined by}
%We defined a new virtual input $\tau$ such that
%$c_0 = \rho S (a_0 - \lambda b_0)$
\begin{equation}
\tau = T \cos{\alpha} + \lambda T \sin{\alpha} - \mu S^\star (a_0 - \lambda b_0 ) T.
\label{eq:tau}
\end{equation}

Since $\gamma^*$ is prescribed \mcedit{by}
%from 
the path, \eqref{eq:nocvx} and \eqref{eq:diff}  form a differential-algebraic system of linear equations and can thus be included as part of a convex program. 

The constraints on the state \mcedit{variables} can be rewritten as 
\begin{gather*}
0 \leq E \leq \overline{V}^2, \quad \underline{a} \leq a \leq \overline{a},
\\
E(s_0) = V_{0}^2
\text{ and }
E(s_f) = V_{f}^2.
\end{gather*}
The thrust constraint \mcedit{$ 0 \leq T \leq \overline{T}$}
%\eqref{eq:thrust_constr} 
cannot be expressed exactly as a function of $\tau$ only. However, assuming that $\alpha \ll 1$, and noting that $\lambda \ll 1$, and $S^\star (a_0 - \lambda b_0 ) \ll 1$ we have $\tau \approx T$ and 
\mcedit{$ 0 \leq T \leq \overline{T}$ is therefore approximately equivalent to}
%the thrust constraint can be approximated by
\[
0\leq \tau \leq \overline{T}. 
%\label{eq:circle}
\]
Finally the minimum thrust criterion in \eqref{eq:obj} can be \mcedit{approximated by}
%turned into 
a convex objective \mcedit{function under these conditions, since $\tau$ is then} 
%if we again assume that $\tau$ is 
a good proxy for the thrust. By the change of differential operator in \eqref{eq:CV} we obtain 
\[
    J = \int_{s_0}^{s_f} \frac{(\tau/\overline{T})^2}{\sqrt{E}} \, \mathrm{d}s.
\]

We can thus state the following convex optimisation problem for a given path (i.e.~given $\gamma$ and $\gamma'$) as follows  
\begin{equation}
\begin{aligned}
& \mathcal{P}_1 : \min_{\substack{\tau,\,E,\,a}} & & \int_{s_0}^{s_f} \frac{(\tau/\overline{T})^2}{\sqrt{E}} \, \mathrm{d}s, 
\nonumber\\
&  \text{ s.t.} & & m a + c(\gamma^{*\prime}) E + d(\gamma^*) = \tau , 
\nonumber\\
& & &    E' = 2 a , 
\nonumber\\
& & &    0 \leq \tau \leq \overline{T},\   \underline{a} \leq a \leq \overline{a}, 
\nonumber\\
& & &   0 \leq E \leq \overline{V}^2,\  E(s_0) = V_{0}^2,\   E(s_f) = V_{f}^2.
\nonumber\\
\end{aligned}
\end{equation}

\mcedit{Solving $\mathcal{P}_1$ yields the optimal velocity}
%Now that we obtained the optimal velocity sequence ($V = \sqrt{E}$) 
along the path and a proxy for the \mcedit{optimal} thrust for sufficiently small values of the angle of attack ($\tau \approx T$).
\mcedit{However, a} 
%we still need to obtain the 
tiltwing angle \mcedit{profile that meets} 
%to realize 
the dynamical constraints and follows the desired path with $\gamma \approx \gamma^*$ \mcedit{must also be computed}.  
\mcedit{To achieve this}
%To do so, 
we use the \mcedit{solution of $\mathcal{P}_1$ to define}
%results from problem $\mathcal{P}_1$ and solve 
a new optimisation problem with variables $\gamma$, $\alpha$, $i_w$,  and $M$ satisfying \mcedit{the constraints}
\begin{gather}
J_w (i_w' a + i_w'' E)  = M,
\label{eq:d2i_w2}
\\
\underline{M} \leq M \leq \overline{M},
\nonumber
\\
i_w  = \alpha + \gamma,
\nonumber
\end{gather}
\small
\begin{gather}
m E \gamma' = \tau \sin{\alpha} + (1-\mu)\frac{1}{2}\rho S (b_1 \alpha + b_0) E -mg \cos{\gamma} 
\nonumber
\\
+\mu\frac{1}{2}\rho S \biggl[ b_1 \arcsin \biggl( \frac{\sqrt{E} \sin{\alpha}}{\sqrt{E + \smash{\frac{2\tau}{\rho A n}}\rule{0pt}{8.5pt}}}\biggr) + b_0\biggr] \Bigl(E + \frac{2\tau}{\rho A n}\Bigr), 
\label{eq:eom_new}
\end{gather}
\normalsize
\mcedit{in which the objective is to minimise the cost function}
%minimizing the following objective 
\begin{equation}
\label{eq:obj2}
    J_{\gamma} = \int_{s_0}^{s_f} \frac{(\gamma - \gamma^*)^2}{\sqrt{E}} \mathrm{d}s. 
\end{equation}
Note that to obtain \eqref{eq:d2i_w2} from \eqref{eq:tiltwing_dyna}, we \mcedit{have} used the result 
$\frac{\mathrm{d^2}}{\mathrm{d}t^2} = V\frac{\mathrm{d}}{\mathrm{d}s} \left( V\frac{\mathrm{d}}{\mathrm{d}s} \right) =a \frac{\mathrm{d}}{\mathrm{d}s} + E \frac{\mathrm{d^2}}{\mathrm{d}s^2}$ inferred from the change of differential operator in \eqref{eq:CV}. To convexify the problem and \mcedit{to}
%also 
support the small angle approximation \mcedit{$\alpha,\alpha_e\ll 1$}
used \mcedit{in the derivation of the EOM},%
\footnote{We refer the reader to \cite{UMich19}, which shows that maintaining a small angle of attack during transition has negligible effect on performance, and this is preferable over operating the wing in dangerous near-stall regimes.} 
%previously, 
we enforce bounds on $\alpha$ by imposing constraints $\underline{\alpha}~\leq~\alpha~\leq~\overline{\alpha}$,
\mcedit{for sufficiently small  $\underline{\alpha}$, $\overline{\alpha}$,} \mbedit{chosen such that the wing is unstalled, and linear relationship with lift coefficient applies}.
This allows \mcedit{the dependence of \eqref{eq:eom_new} on $\alpha$} to be linearised %around $\alpha=0$ as follows
\small
\begin{multline}
m E \gamma' = \underbrace{\biggl(\tau + (1\!-\!\mu)\frac{1}{2}\rho S b_1 E   +\mu\frac{1}{2}\rho S b_1 \sqrt{E^2 + \tfrac{2\tau E}{\rho A n}} \biggr)}_{p(E,\tau)}\alpha \\
+ \underbrace{(1 - \mu)\frac{1}{2}\rho S b_0 E +\mu\frac{1}{2}\rho S b_0 \Bigl(E + \frac{2\tau}{\rho A n}\Bigr)}_{q(E, \tau)} -mg \cos{\gamma}.
\label{eq:eom_small}
\end{multline}
\normalsize
\mcedit{Although $p(E,\tau)$ and $q(E,\tau)$ are determined by the solution of $\mathcal{P}_1$, \eqref{eq:eom_small} is nonconvex since} the last term is nonlinear in the \mcedit{decision} variable $\gamma$. However, 
\mcedit{\eqref{eq:eom_small} can be omitted  from the problem if the objective function \eqref{eq:obj2} is augmented}
%\mcedit{if a solution exists satisfying $\gamma = \gamma^\ast$, then the}
%we can get rid of it by noting that 
%objective in \eqref{eq:obj2} is equivalent to the augmented objective
 \begin{equation}
\label{eq:obj3}
       J^a_{\gamma} = \int_{s_0}^{s_f} \frac{(\gamma - \gamma^*)^2}{\sqrt{E}} + 
\frac{(p \alpha + q - m E \gamma' -m g \cos{\gamma^*})^2}{(mg)^2\sqrt{E}} \mathrm{d}s .
%\frac{( \cos{\gamma} -\cos{\gamma^*})^2}{\sqrt{E}} \mathrm{d}s. 
\end{equation}
This \mcedit{objective attempts to minimise}
%choice allows to minimize 
the error between the \mdsedit{guessed} \mcedit{and actual}
%flight path angle and current 
flight path angles \mcedit{while enforcing 
%while allowing to include 
\eqref{eq:eom_small} in the sense that, if a feasible solution exists satisfying both $\gamma = \gamma^\ast$ and \eqref{eq:eom_small}, then the optimal solution must satisfy \eqref{eq:eom_small}}. 

We thus state the following optimisation problem, with $\tau$, $E$ and $a$ prescribed from problem $\mathcal{P}_1$
\begin{equation}
\begin{aligned}
& \mathcal{P}_2 : \min_{\substack{\alpha,\,\gamma,\,i_w,\, M}} & & \int_{s_0}^{s_f} \frac{(\gamma - \gamma^*)^2}{\sqrt{E}} \mathrm{d}s \\&&&+ \int_{s_0}^{s_f}\frac{(p \alpha + q - m E \gamma' -m g \cos{\gamma^*})^2}{(mg)^2\sqrt{E}} \mathrm{d}s
\nonumber\\
&  \text{ s.t.} & & J_w (i_w' a + i_w'' E)  = M,\ i_w(s_0) = i_{0}, 
\nonumber\\
& & &   i_w'(s_0) \sqrt{E(s_0)}  = \Omega_{0},
\nonumber\\
& & &    i_w  = \alpha + \gamma, 
\nonumber\\
& & &    \underline{M} \leq M \leq \overline{M},\   \underline{\alpha} \leq \alpha \leq \overline{\alpha}
\nonumber\\
& & &   \underline{\gamma} \leq \gamma \leq \overline{\gamma},\  \underline{i_w}  \leq i_w \leq \overline{i_w}. 
\nonumber\\
\end{aligned}
\end{equation}
%Note that inequality constraints on $\gamma$ and $i_w$ \mcedit{have been} added for improved numerical stability. 

It should be emphasized that $\mathcal{P}_2$ may admit solutions that allow $\gamma$ to divert 
%freely 
from the \mdsedit{guessed} flight path angle $\gamma^*$. This occurs for example when the bounds on angle of attack do not allow 
%to achieve 
the desired path \mcedit{to be followed (see e.g.~Section \ref{sec:results}, scenario 2)}.  In \mcedit{this case},
%such event, 
problems $\mathcal{P}_1$ and $\mathcal{P}_2$ can be 
%\mcedit{solved} successively, with $\gamma^\ast$ redefined 
rerun sequentially after reinitializing with the newly obtained path ($\gamma^* \gets \gamma$). This process can be repeated iteratively until $\gamma$ is sufficiently close to $\gamma^*$. 
%\mcedit{We  no theoretical guarantee The resulting iteration 
Although we \mcedit{do not} provide theoretical guarantees of convergence \mcedit{for this} case, \mcedit{we note} 
%observe empirically 
that the objective value \mcedit{has decreased by orders of magnitude}
%significantly in magnitude 
after several iterations, as illustrated in Section \ref{sec:results}. 

%It is worthwhile mentioning that 
Only one EOM \mcedit{(\ref{eq:eom8})} was needed to constrain the variables of problem $\mathcal{P}_2$ since the EOM \mcedit{(\ref{eq:eom7})} can be enforced \textit{a posteriori} by choosing $T$ so that equation \eqref{eq:tau} is satisfied
\[
%    \label{eq:T}
    T = \frac{\tau}{ \cos{\alpha} + \lambda  \sin{\alpha} - \mu S^\star (a_0 - \lambda b_0 )}.  
\]
This guarantees that the solution satisfies both EOM within the required accuracy,
%\footnote{Provided $\alpha$ is kept below a certain threshold, which can be enforced explicitly with the constraint on the angle of attack in problem $\mathcal{P}_2$.}, 
since problem $\mathcal{P}_1$ includes a linear combination of \mcedit{(\ref{eq:eom7}) and (\ref{eq:eom8}), and problem $\mathcal{P}_2$ enforces (\ref{eq:eom8})}.

\mcedit{Given an optimal set of} inputs and states as functions of the \mcedit{independent variable} $s$, the \mcedit{final} step is to convert back to time domain. This can be done by \mcedit{reversing} the change of differential operator in equation \eqref{eq:CV} and integrating%, which allows us to define time $t$ as a parametric integral
\[
%\label{eq:back2time}
t(\xi) = \int_{s_0}^\xi \frac{\mathrm{d}s }{V(s)}.
\]
%
%A time $t(\xi)$ is thus associated with each point $\xi$ on the path.

\section{Discrete convex optimisation}
\label{sec:discrete}
\mcedit{The decision variables in} $\mathcal{P}_1$ and $\mathcal{P}_2$ \mcedit{are functions defined on the interval $[s_0,s_f]$}.
%are continuous time, infinite dimensional convex optimization problems. 
To obtain \mcedit{computationally tractable problems}, we 
\mcedit{consider}
%discretise the problem, considering 
$N+1$ discretisation points $\{s_0, s_1,  \ldots, s_{N}\}$ of the path, with spacing $\delta_k = s_{k+1} - s_{k}$, $k = 0, \ldots, N-1$ ($N$ steps). 
%\mcedit{chosen} so as to refine the mesh in portions of the path where the \mcedit{aircraft} velocity is \mcedit{expected} to be \mcedit{lowest}. \mcedit{This refinement compensates for the fact that}, according to \eqref{eq:CV}, \mcedit{sections} of the path with low velocities yield larger time steps for given \mcedit{spatial} step sizes. 
%
The notation $\{u_0, \ldots, u_{N}\}$ is used for the sequence of the discrete values of a continuous variable $u$ evaluated at the discretisation points of the mesh, where $u_k = u(s_{k})$, $\forall k \in \{0, \ldots, N\}$.

Assuming a path $s_k \rightarrow (x_k, z_k)$, the prescribed flight path angle and rate are discretised as follows
\begin{align}
\gamma_k^* &= \arctan \Bigl(-\frac{z_{k+1}- z_{k}}{x_{k+1}- x_{k}}\Bigr), \quad k \in \{0, \ldots, N-1\},
\label{eq:gamma}
\\
\gamma_k^{*\prime} &= \begin{cases} (\gamma_{k+1}^*- \gamma_{k}^*)/{\delta_k}, &  k \in \{0, \ldots, N-2\}, \\
\gamma_{N-2}^{*\prime}, & k = N-1.
\end{cases}
\label{eq:dgamma}
\end{align}
The resulting discrete time versions of $\mathcal{P}_1$ and $\mathcal{P}_2$ are
\begin{equation}
\begin{aligned}
& \mathcal{P}_1^\dagger : \min_{\substack{\tau,\,E,\,a}} & & \sum_{k=0}^{N-1} \frac{(\tau_k/\overline{T})^2}{\sqrt{E_k}} \, \delta_k,
\nonumber\\
&  \text{ s.t.} & & m a_k + c(\gamma_k^{*\prime}) E_k + d(\gamma_k^*) = \tau_k , 
\nonumber\\
& & &    E_{k+1} = E_{k} +  2 a_{k} \delta_{k}, 
\nonumber\\
& & &    0 \leq \tau_k \leq \overline{T}, \quad  \underline{a} \leq a_k \leq \overline{a}, 
\nonumber\\
& & &   0 \leq E_k \leq \overline{V}^2, \quad  E_0 = V_{0}^2, \quad   E_{N} = V_{f}^2,
\nonumber\\
\end{aligned}
\end{equation}
\begin{equation}
\begin{aligned}
& \mathcal{P}_2^\dagger : \min_{\substack{\alpha,\,\gamma,\, \theta,\\i_w,\, \zeta,\, M}} & & \sum_{k=0}^{N-1} \frac{(\gamma_k - \gamma_k^*)^2}{\sqrt{E_k}}\delta_k \\ & & & + \sum_{k=0}^{N-1} \frac{(p \alpha_k + q - m E_k \Psi_k -m g \cos{\gamma_k^*})^2}{(mg)^2\sqrt{E_k}} \delta_k, 
\nonumber\\
& \text{s.t.} & &  \gamma_{k+1} =  \gamma_{k} + \Psi_k \delta_k , 
\nonumber\\
& & &    i_{w,k}  = \alpha_k + \gamma_k, 
\nonumber\\
& & & i_{w,k+1}   = i_{w,k} + \zeta_k \delta_k, \quad i_{w,0} = i_{0}, 
\nonumber\\
& & & \zeta_{k+1}  = \zeta_{k}\Bigl(1-\frac{a_k \delta_k}{E_k}\Bigr) + \frac{M_k \delta_k}{J_w {E_k}}, 
\nonumber \\
& & & \zeta_0 \sqrt{E}_0  = \Omega_{0},
\nonumber \\
& & &    \underline{M} \leq M_k \leq \overline{M}, \quad \underline{\alpha} \leq \alpha_k \leq \overline{\alpha}, 
\nonumber\\
& & &   \underline{\gamma} \leq \gamma_k \leq \overline{\gamma}, \quad   \underline{i_w}  \leq i_{w,k} \leq \overline{i_w}.
\nonumber\\
\end{aligned}
\end{equation}
%Here $\Psi$ is introduced as a dummy variable to decouple the derivative for convenience.

Once 
%problems 
$\mathcal{P}_1^\dagger$ and $\mathcal{P}_2^\dagger$ have been solved, we check whether $|\gamma_k^*-\gamma_k| \leq \epsilon$ $\forall k \in \{0, \ldots, N\}$, where $\epsilon$ is a specified tolerance. If this condition is not met, the problem is reinitialized with the updated flight path angle and rate $\gamma^*_k \gets \gamma_k$, $\gamma^{*\prime}_k \gets \gamma_k'$.
and \mcedit{$\mathcal{P}_1^\dagger$ and $\mathcal{P}_2^\dagger$ are solved again}.
% until convergence. 
When \mcedit{the solution tolerance is met}
%the test is passed 
(or 
%after a given 
\mcedit{the} maximum number of iterations \mcedit{is exceeded}) the problem is considered solved and the 
%real 
input and state \mcedit{variables} are reconstructed 
using
%as follows
\begin{align*}
    T_k &= \frac{\tau_k}{ \cos{\alpha_k} + \lambda  \sin{\alpha_k} - \mu S^\star (a_0 - \lambda b_0 )}, \quad V_k &= \sqrt{E_k} ,
\end{align*}
\mcedit{and} 
%along with 
the time $t_k$ associated with each discretisation point 
\mcedit{is computed, allowing solutions to be expressed}
%to express results 
as time series
\[
\label{eq:back2time}
t_k = \sum_{j=0}^{k-1} \frac{\delta_j}{V_j}.
\]
The procedure is summarised in Algorithm 1. 

%%%%%%%%%%%%%%%%%%%%%%%%%%%%%%%%%%%%%%%%%%%%%%%%%%%%%%%%%%%%%%%%%%%%%%%%%%%%%%%%%%%%%%%%%%%%%%%%%%%%%%%%%%%%%%%%%%%%%%%%%%%%%%%%%%%%%%%%%%%
% Martin: The way I initialised the algorithm is a bit far-fetched. I use equations (16-17) that are given for prescribed variables in order to initialise current variables, then assign the current variables to prescribed variables inside the while loop. I was thinking of using a "while True" with a conditional break statement inside it, but you might think of something better to present this ?
% Mark: To be clearer we could say explicitly that gamma and gamma' are set equal to the gamma^* and gamma'^* computed using (16) and (17) - see the updated Algorithm 1
%%%%%%%%%%%%%%%%%%%%%%%%%%%%%%%%%%%%%%%%%%%%%%%%%%%%%%%%%%%%%%%%%%%%%%%%%%%%%%%%%%%%%%%%%%%%%%%%%%%%%%%%%%%%%%%%%%%%%%%%%%%%%%%%%%%%%%%%%%%

\begin{algorithm}
\caption{Convex trajectory optimisation}
\label{algo}
%\small
    %\SetAlgoLined
    \mcedit{Compute $\gamma^*$, $\gamma^{*\prime}$ using \eqref{eq:gamma}, \eqref{eq:dgamma}} and \mcedit{initialise: $\gamma\gets\gamma^*$, $\gamma'\gets\gamma^{*\prime}$,} $\gamma^* \gets \infty$, $j \gets 0$\\
    \While{$\displaystyle\max_{k \in \{0, \ldots, N\}}|\gamma_k^*-\gamma_k| > \epsilon$ $\&$ $j < $ MaxIters}{
    $\gamma^{*} \gets \gamma$, $\gamma^{*\prime} \gets \gamma'$\\
    Solve problem $\mathcal{P}_1^\dagger$ and problem $\mathcal{P}_2^\dagger$ 
%successively. 
\\
    $j \gets j+1$
    }
   
    \For{$k\gets0$ \KwTo $N$}{
    $T_k \gets \frac{\tau_k}{ \cos{\alpha_k} + \lambda  \sin{\alpha_k} - \mu S^\star (a_0 - \lambda b_0 )}$ \\
    $V_k \gets \sqrt{E_k}$\\
    $t_k \gets \sum_{j=0}^{k-1} \frac{\delta_j}{V_j}$
    }
\end{algorithm}

\section{Results}
\label{sec:results}
We consider a case study \mcedit{based on}
%inspired from 
the Airbus A$^3$ Vahana. The \mcedit{aircraft}
%vehicle 
parameters are reported in Table \ref{tab:param}.
We run Algorithm 1 \mcedit{using}
%and use 
the convex programming software package CVX \cite{cvx} with \mcedit{the solver} Mosek \cite{mosek} to compute the optimal trajectory for 3 different transition manoeuvres, with boundary conditions given in Table \ref{tab:BC}. 

\begin{figure}
         \centering
         \includegraphics[width=0.5\textwidth, trim={0cm 0.3cm 0cm 0.5cm}, clip]{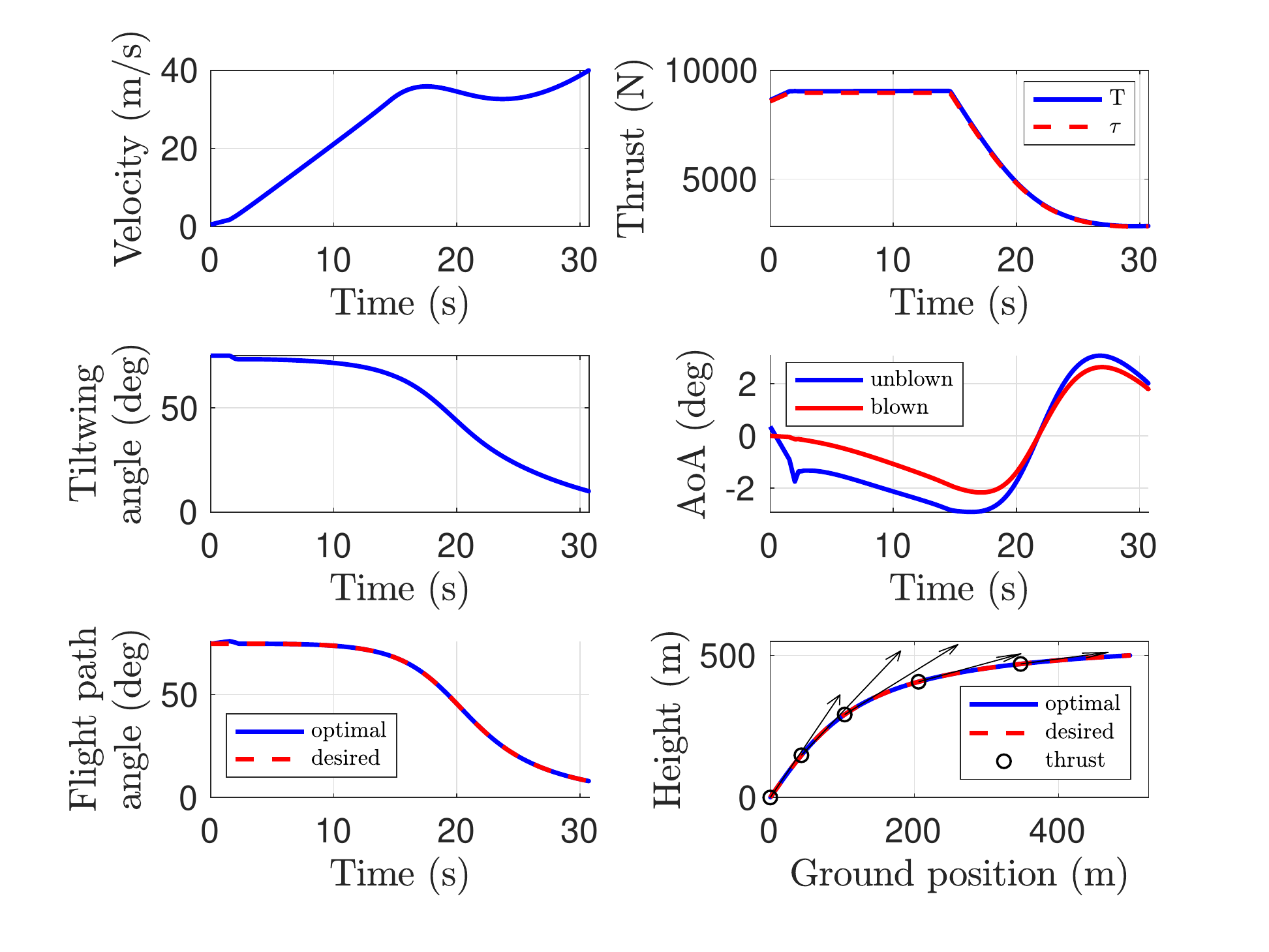}
         \caption{Forward transition (scenario 1)}
         \label{fig:traj1}
\end{figure}

The first \mcedit{scenario} is a forward smooth transition where the flight path angle varies slowly from $75^\circ$ to 
\mcedit{a value close to zero},
%small values, 
as illustrated in Figure \ref{fig:traj1}. As the \mcedit{aircraft transitions}
%vehicle switches 
from powered lift
%flight 
to cruise, the velocity magnitude increases and the thrust decreases, illustrating the change in lift \mcedit{generation} from propellers to wing. The angle of attack remains small over the transition as it occurs smoothly. 

\begin{figure}
         \centering
         \includegraphics[width=0.4\textwidth, trim={0cm 0cm 0cm 0.9cm}, clip]{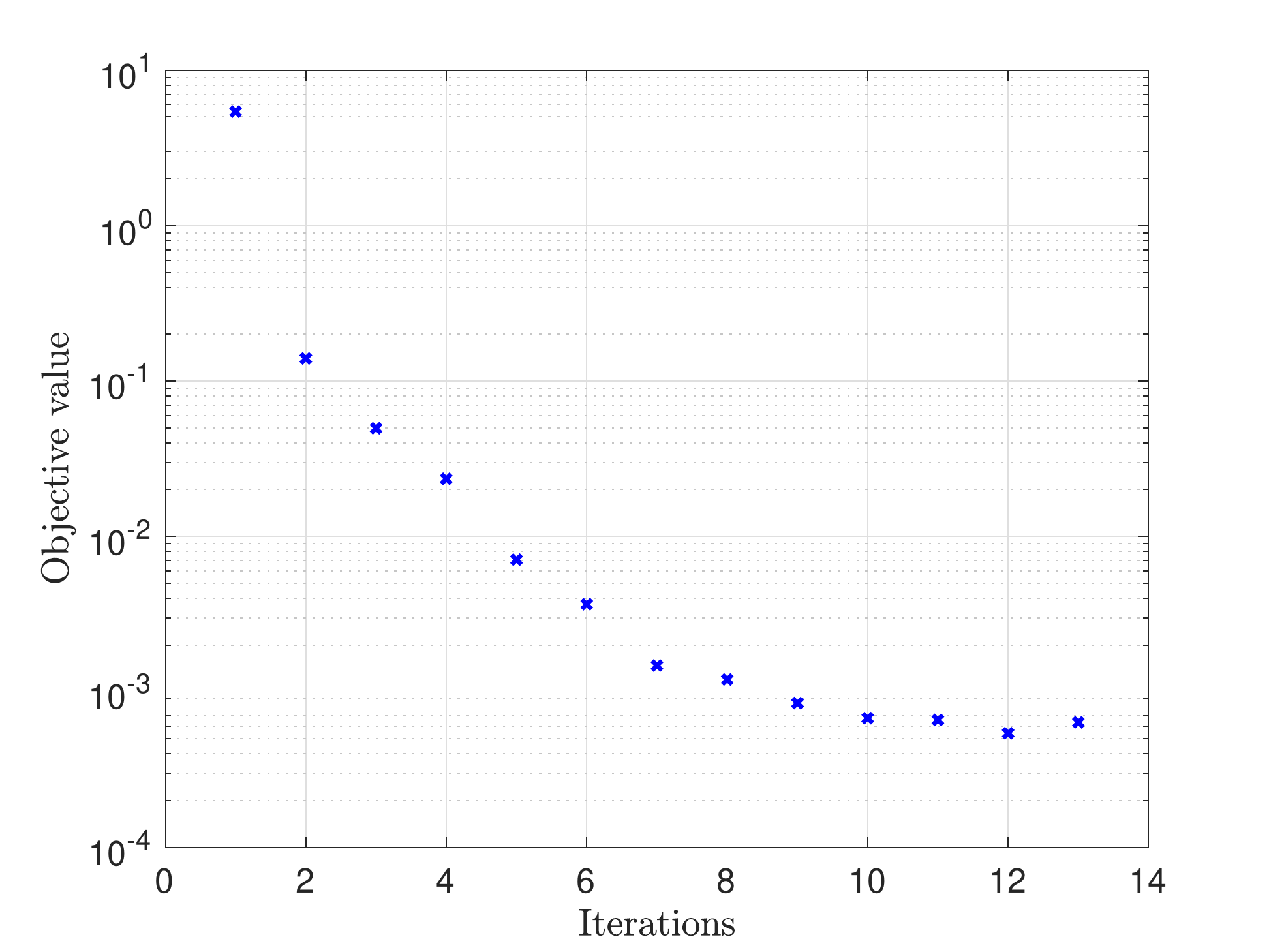}
         \caption{Convergence of objective value for $\mathcal{P}_2^\dagger$ (scenario 2)}
         \label{fig:obj}
\end{figure}

\begin{figure}
         \centering
         \includegraphics[width=0.5\textwidth, trim={0cm 0cm 0cm 0cm}, clip]{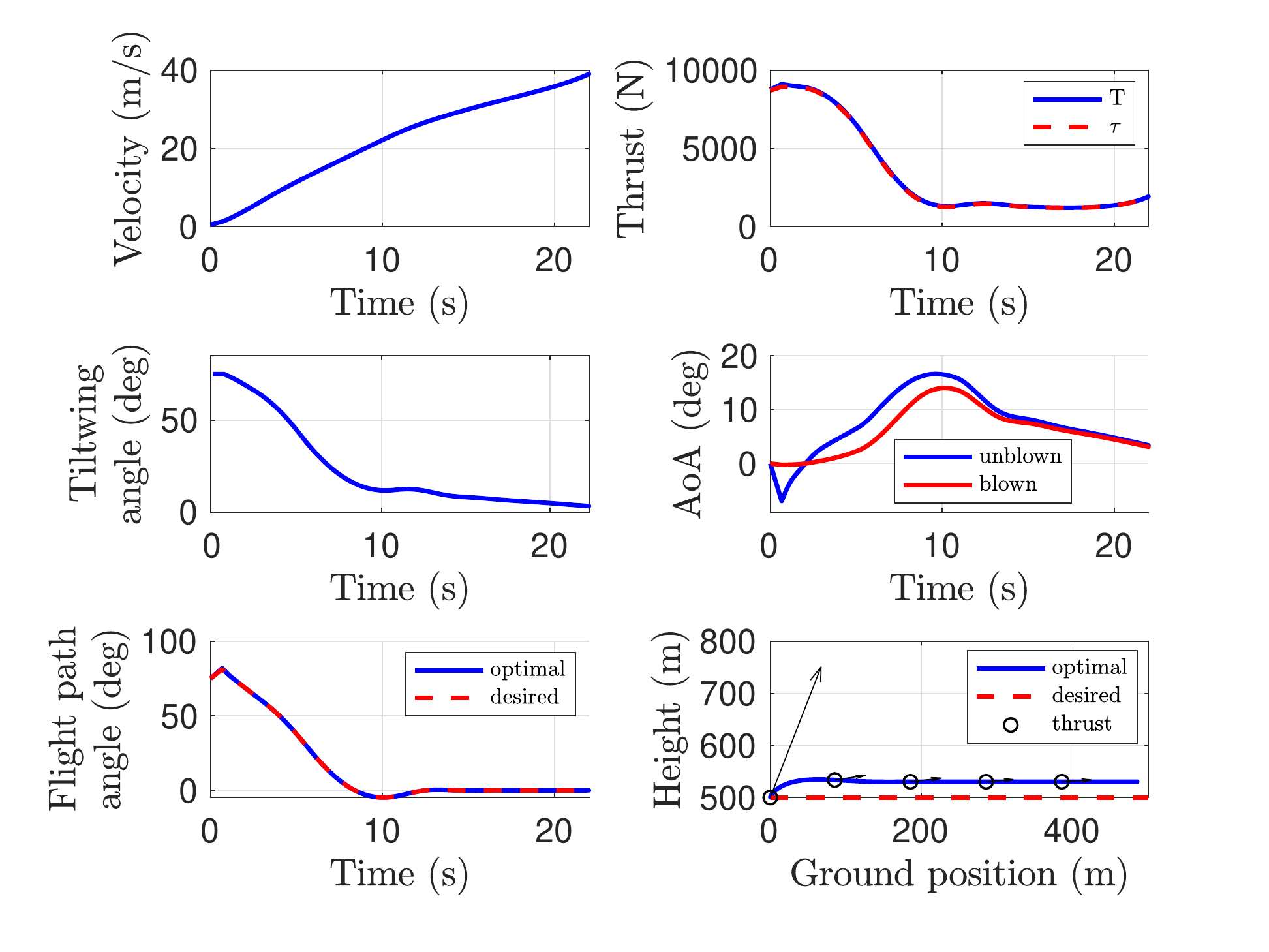}
         \caption{Forward transition (scenario 2)}
         \label{fig:traj2}
\end{figure}

The second \mcedit{scenario} is a constant altitude forward transition. This manoeuvre is more abrupt and would require a zero flight path angle throughout, thus exceeding the bounds on $\alpha$ at the beginning when the wing is tilted vertically. Solving problem $\mathcal{P}_2^\dagger$ with 
\mcedit{constraints}
%such bounds 
on $\alpha$ 
\mcedit{initially yields}
%first yielded 
an optimal $\gamma$ profile 
\mcedit{that differs}
%different 
from the \mcedit{anticipated}
%guessed 
zero flight path angle. 
\mcedit{Reinitializing}
%This required to reinitialize 
the problem with the \mcedit{computed}
%obtained 
flight path angle ($\gamma^* \gets \gamma$) and 
\mcedit{iterating}
%iterate 
until a satisfactory match between $\gamma$ and $\gamma^*$ 
\mcedit{is}
%was 
achieved 
\mcedit{results in convergence illustrated empirically in Figure \ref{fig:obj}}. 
%After a few iterations, 
\mcedit{The solution thus obtained is shown} 
%collected 
in Figure \ref{fig:traj2}. Since the 
%achieved 
flight path angle is not identically zero (as would be required for a constant altitude transition) the manoeuvre is characterised by an initial increase in altitude and overshoot of about 30\,m. This is 
\mcedit{a consequence of}
%the price to pay for 
limiting the angle of attack. More abrupt transitions with %no 
\mcedit{smaller altitude variations} are possible, but this \mcedit{requires stalling}
%would require to stall 
the wing and \mcedit{operating} in flight envelopes \mcedit{for which our approximations are not valid}. 
%(and hence where the present solver becomes obsolete). 

\begin{figure}
         \centering
         \includegraphics[width=0.5\textwidth, trim={0cm 0cm 0cm 0cm}, clip]{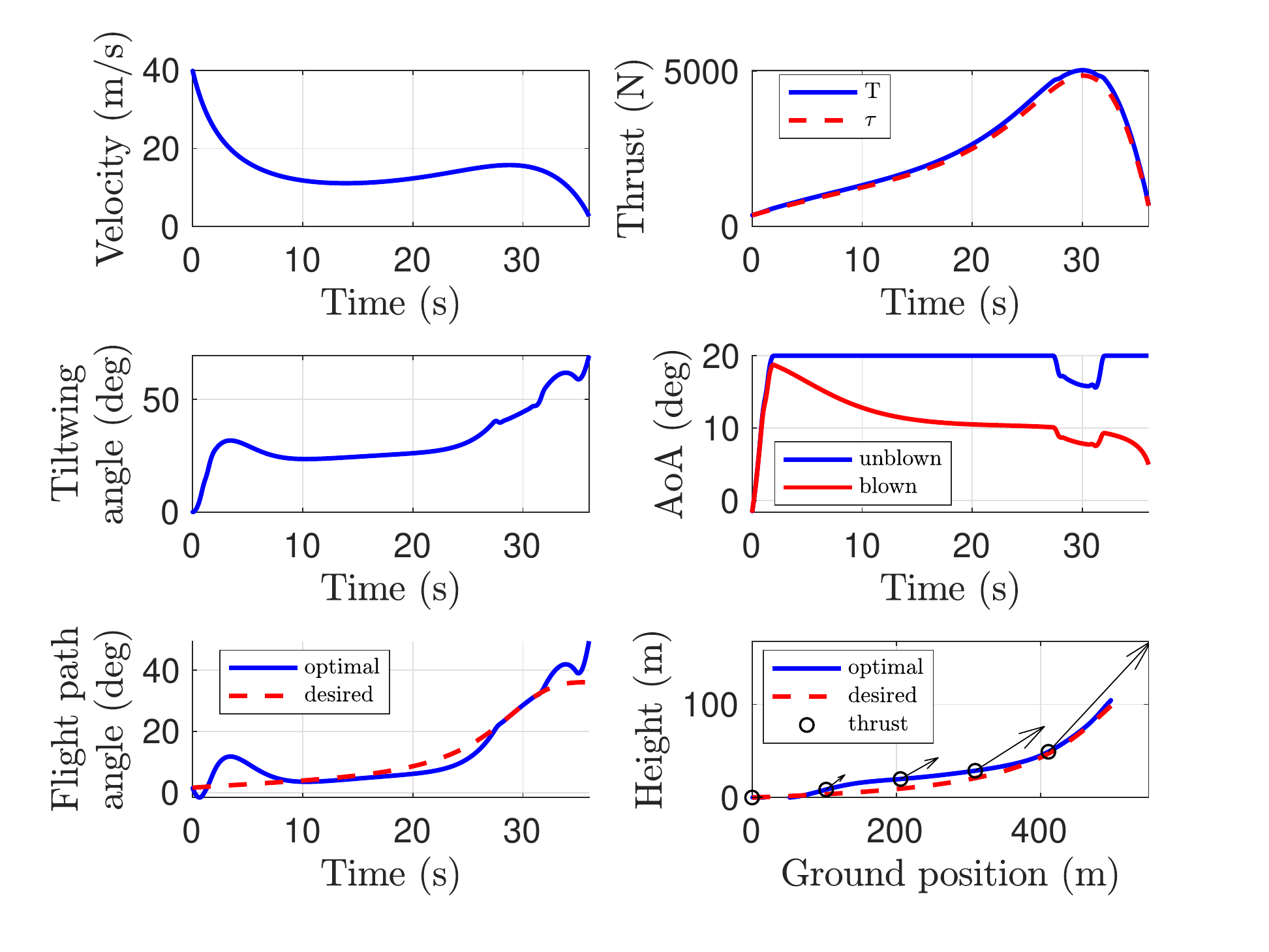}
         \caption{Backward transition (scenario 3)}
         \label{fig:traj3}
\end{figure}

The third \mcedit{scenario} is a backward transition with an increase in altitude (Figure \ref{fig:traj3}). 
\mcedit{This}
%It 
is characterised by \mcedit{an initial}
%a 
decrease in velocity magnitude and increase in thrust.
%throughout. 
The strict bounds on angle of attack \mcedit{necessitate an increase in altitude of about 100\,m}.
%require to consider manoeuvres with large gains in altitude as illustrated. 
A backward transition at constant altitude would require \mcedit{stalling} 
%to stall 
the wing, which is prohibited in the present formulation, illustrating a limitation of our approach. 
%\mcedit{In order}
%Note that 
To achieve the backward transition, 
%the deployment of 
a high-drag device or flaps \mcedit{are needed to provide braking forces}.
%is required to brake. 
This can be modelled by adding a constant term $\delta$ to $c(\gamma^{*\prime})$  in problem $\mathcal{P}_1^\dagger$ when simulating the backward transition.

\begin{figure}
         \centering
         \includegraphics[width=0.4\textwidth, trim={0cm 0cm 0cm 0cm}, clip]{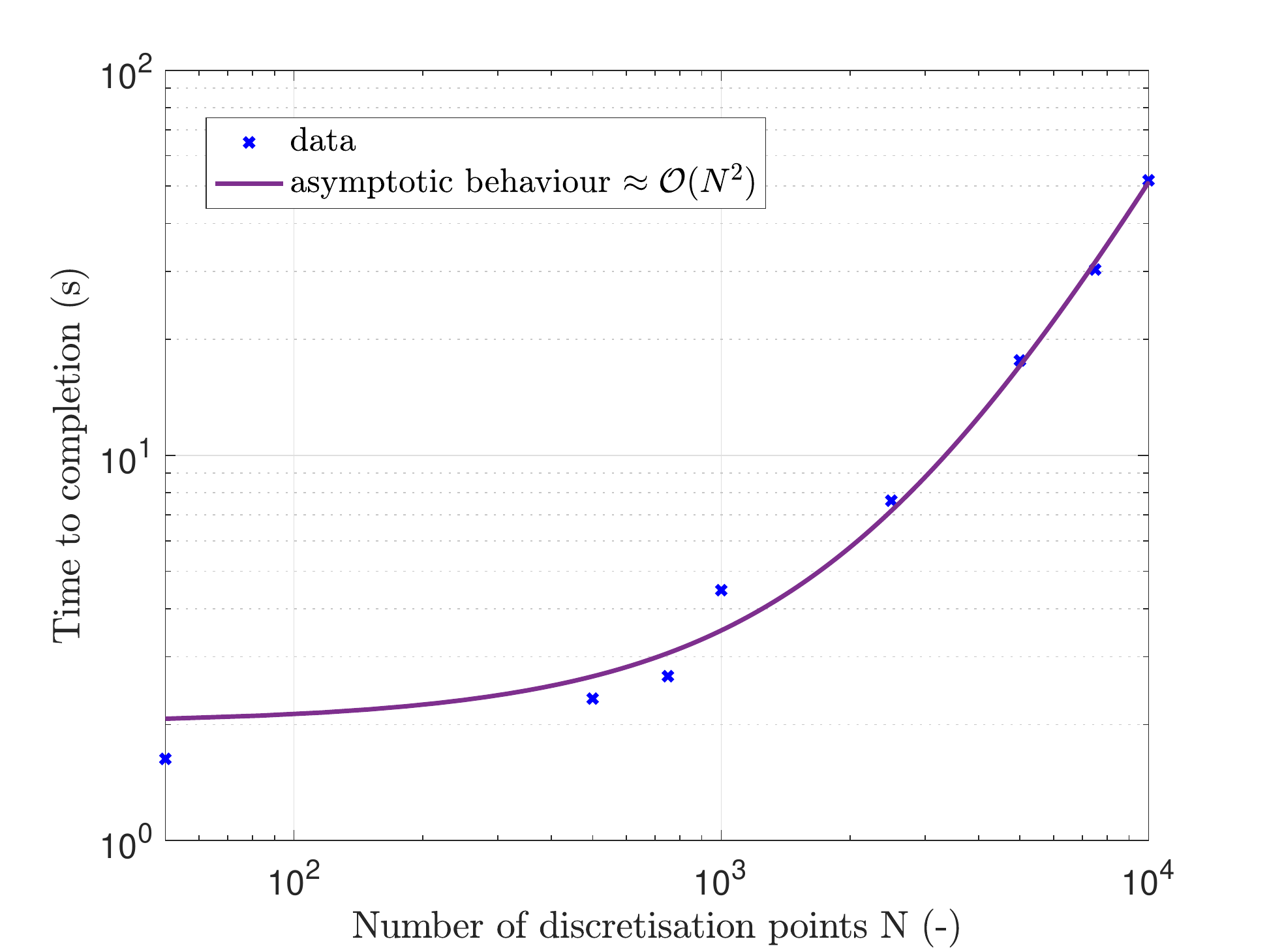}
         \caption{Time to completion as a function of problem size}
         \label{fig:time}
\end{figure}

The time-complexity of the
\mcedit{proposed algorithm in scenario~1 is shown in Figure \ref{fig:time}}.
%present solver. 
The time to completion \mcedit{as a function of  problem size (number of discretisation points) is asymptotically quadratic, implying excellent scalability}. In particular, for $N<1500$, the time to completion is of the order of seconds. A real-time solution is thus possible 
\mcedit{using modest computational resources}. 
\mcedit{This is in stark contrast to the complexity of generic NLP approaches \cite{leo}.}
%on modern embedded systems. 
\mdsedit{Moreover, as shown in \cite{me3}, the time complexity for this type of problem can be reduced significantly with first-order solvers, paving the way for fast real-time implementations}. 

\begin{comment}
\subsection{Forward smooth manoeuvre}
To be continued...
\subsection{Forward constant altitude manoeuvre}
To be continued...
\subsection{Backward manoeuvre}
To be continued...
\end{comment}

\begin{table}[ht]
\centering
\begin{tabular}{llll}
\hline
\textbf{Parameter} & \textbf{Symbol} & \textbf{Value} & \textbf{Units} \\ \hline
Mass     &    $m$    &   $752.2$    &    \si{kg}   \\ \hline
Gravity acceleration     &    $g$    &   $9.81$    &    \si{m.s^{-2}}   \\ \hline
Wing area &   $S$     &    $8.93$   &    \si{m^2}   \\ \hline
Disk area &   $A$     &    $2.83$   &    \si{m^2} \\ \hline
Blown-unblown ratio &   $\mu$     &    $0.73$   &    \si{-}   \\ \hline
Wing inertia &   $J_w$     &    $1100$   &    \si{kg.m^2}   \\ \hline
%Induced speed coefficient &   $k_w$     &    $2$   &    \si{-} \\ \hline
Density of air &   $\rho$     &    $1.225$   &    \si{kg.m^{-3}} \\ \hline
%\makecell{Specific heat of air \\(at constant pressure)}    &    $c_p$    &   $1000$    &    \si{J.K^{-1}.kg^{-1}}   \\ \hline
{Lift coefficients} &   $b_0, b_1$   &    $0.43, 0.11$   &    \si{-},  \si{deg^{-1}}  \\ \hline
%&   $b_1$     &    $0.11$   &    \si{deg^{-1}}  \\ \hline
 {Drag coefficients} &   $a_0, a_1$   &   $0.029, 0.004$   &    \si{-}, \si{deg^{-1}}  \\ \hline
%&   $a_1$   &   $0.004$   &    \si{deg^{-1}} \\ \hline
%&   $a_2$   &   $5.3\mathrm{e}{-4}$   &    \si{deg^{-2}}\\ \hline
Maximum thrust &   $\overline{T}$     &    $8855$   &    \si{N} \\ \hline
Angle of attack range &   $\left[\underline{\alpha}, \overline{\alpha}\right]$     &    $\left[-20, 20\right]$   &    \si{deg}   \\ \hline
Flight path angle range &   $\left[\underline{\gamma}, \overline{\gamma}\right]$     &    $\left[-90, 90\right]$   &    \si{deg}   \\ \hline
Tiltwing angle range &   $\left[\underline{i}_w, \overline{i}_w\right]$     &    $\left[0, 100\right]$   &    \si{deg}   \\ \hline
Acceleration range &   $\left[\underline{a}, \overline{a}\right]$     &    $\left[-0.3g, 0.3g\right]$   &    \si{m . s^{-2}}   \\ \hline
Velocity range &   $\left[\underline{V}, \overline{V}\right]$     &    $\left[0, 40\right]$   &    \si{m/s}   \\ \hline
Momentum range &   $\left[\underline{M}, \overline{M}\right]$     &    $\left[-50; 50\right]$   &    \si{N.m}   \\ \hline
%Velocity boundary conditions &   $\left\{V_0; V_f\right\}$     &    $\left\{1, 40\right\}$   &    \si{m/s}   \\ \hline
Number of propellers &   $n$     &    $4$   &    \si{-} \\ \hline
Discretisation points &   $N$     &    $1500$   &    \si{-} \\ \hline 
\end{tabular}
\vspace{1mm}\caption{Model parameters \mcedit{derived} from A$^3$ Vahana}
\label{tab:param}
\vspace{-3mm}
\end{table}

\begin{table}[ht]
\centering
\begin{tabular}{llll}
\hline
\textbf{Parameter} & \textbf{Symbol} & \textbf{Value} & \textbf{Units} \\ \hline
\multicolumn{4}{c}{{\cellcolor[rgb]{0.753,0.753,0.753}}\textbf{Forward transition}}    \\
Velocity  &   $\left\{V_0; V_f\right\}$     &    $\left\{0.5; 40\right\}$   &    \si{m/s}   \\ \hline
Tiltwing angle  &   $i_0$     &    $75$   &    \si{deg}   \\ \hline
Tiltwing angle rate  &   $\Omega_0$     &    $0$   &    \si{deg/s}   \\ \hline
Flight path angle &   $\gamma_0$     &    $75$   &    \si{deg}   \\ \hline
\multicolumn{4}{c}{{\cellcolor[rgb]{0.753,0.753,0.753}}\textbf{Backward transition}}    \\
Velocity  &   $\left\{V_0; V_f\right\}$     &    $\left\{40; 0.1\right\}$   &    \si{m/s}   \\ \hline
Tiltwing angle  &   $\left\{i_0; i_f\right\}$     &    $\left\{0; 75\right\}$    &    \si{deg}   \\ \hline
Tiltwing angle rate  &   $\Omega_0$     &    $0$   &    \si{deg/s}   \\ \hline
Flight path angle &   $\gamma_0$     &    $1.6$   &    \si{deg}   \\ \hline
\end{tabular}
\vspace{1mm}\caption{Boundary conditions for transitions}
\label{tab:BC}
\vspace{-3mm}
\end{table}

\section{Conclusions}
\label{sec:conclusion}
 A convex programming formulation of the minimum thrust transition problem for a \mcedit{tiltwing} VTOL aircraft with propeller-wing interaction was proposed and solved for various transition scenarios. The %present 
 approach can compute an optimal trajectory within seconds for large numbers of discretisation points and is particularly efficient at computing smooth transitions. A potential limitation of the approach is its reliance on small angles of attack, which may restrict the range of achievable manoeuvres. In particular, more abrupt manoeuvres with larger angles of attack are prohibited. \mcedit{Another drawback is that, for the case in which the initial guess of flight path angle is infeasible, the evidence for the convergence of the proposed iteration is empirical rather than theoretical}.
 %convergence is not guaranteed when 
 %the initial guess of flight path angle is
 %not feasible, requiring to run the problem iteratively. 
Future work will address these two issues using \mcedit{techniques} %concepts 
from robust MPC theory. 
A further extension to this work 
%would include developing 
\mcedit{is to develop a first-order solver for online trajectory generation}.

\bibliography{biblio} 

\begin{thebibliography}{10}

\bibitem{easa}
McKinsey, ``Study on the societal acceptance of {U}rban {A}ir {M}obility in
  {E}urope,'' {\em European Union Aviation Safety Agency (EASA)}, 2021.

\bibitem{roland-berger}
M.~Hader, S.~Baur, S.~Kopera, T.~Schönberg, and J.-P. Hasenberg, ``Urban air
  mobility, {USD} 90 billion of potential: how to capture a share of the
  passenger drone market,'' {\em Roland Berger}, 2020.

\bibitem{kuhn1959semiempirical}
R.~E. Kuhn, ``Semiempirical procedure for estimating lift and drag
  characteristics of propeller-wing-flap configurations for vertical-and
  short-take-off-and-landing airplanes,'' {\em {NASA} Memorandum}, 1959.

\bibitem{hargraves1961analytical}
C.~R. Hargraves, ``An analytical study of the longitudinal dynamics of a
  tilt-wing {VTOL},'' tech. rep., Princeton University, 1961.

\bibitem{mccormick1967aerodynamics}
B.~W. McCormick, {\em Aerodynamics of {V}/{STOL} flight}.
\newblock Academic Press, London, 1967.

\bibitem{UMich19}
S.~S. Chauhan and J.~R. Martins, ``Tilt-wing e{VTOL} takeoff trajectory
  optimization,'' {\em Journal of Aircraft}, pp.~1--20, 2019.

\bibitem{rothhaar2014nasa}
P.~M. Rothhaar, P.~C. Murphy, B.~J. Bacon, I.~M. Gregory, J.~A. Grauer, R.~C.
  Busan, and M.~A. Croom, ``{NASA} langley distributed propulsion {VTOL}
  tiltwing aircraft testing, modeling, simulation, control, and flight test
  development,'' in {\em 14th AIAA aviation technology, integration, and
  operations conference}, p.~2999, 2014.

\bibitem{higgins2020aeroacoustic}
R.~J. Higgins, G.~N. Barakos, S.~Shahpar, and I.~Tristanto, ``An aeroacoustic
  investigation of a tiltwing e{VTOL} concept aircraft,'' in {\em AIAA AVIATION
  2020 FORUM}, p.~2684, 2020.

\bibitem{me}
R.~Chiappinelli, M.~Cohen, M.~Doff-Sotta, M.~Nahon, J.~R. Forbes, and
  J.~Apkarian, ``Modeling and control of a passively-coupled tilt-rotor
  vertical takeoff and landing aircraft,'' in {\em 2019 International
  Conference on Robotics and Automation (ICRA)}, pp.~4141--4147, IEEE, 2019.

\bibitem{Iowa19}
P.~Pradeep and P.~Wei, ``Energy optimal speed profile for arrival of tandem
  tilt-wing e{VTOL} aircraft with {RTA} constraint,'' in {\em IEEE CSAA
  Guidance, Navigation and Control Conference}, 2018.

\bibitem{leo}
L.~Panish and M.~Bacic, ``Transition trajectory optimization for a tiltwing
  {VTOL} aircraft with leading-edge fluid injection active flow control,'' {\em
  {AIAA} Scitech 2022 San Diego}, 2022.

\bibitem{bobrow}
J.~E. Bobrow, S.~Dubowsky, and J.~S. Gibson, ``Time-optimal control of robotic
  manipulators along specified paths,'' {\em The international journal of
  robotics research}, vol.~4, no.~3, pp.~3--17, 1985.

\bibitem{me2}
M.~Doff-Sotta, M.~Cannon, and M.~Bacic, ``Optimal energy management for hybrid
  electric aircraft,'' {\em IFAC-PapersOnLine}, vol.~53, no.~2, pp.~6043--6049,
  2020.

\bibitem{me3}
M.~Doff-Sotta, M.~Cannon, and M.~Bacic, ``Predictive energy management for
  hybrid electric aircraft propulsion systems,'' {\em {IEEE} Transactions on
  Control Systems Technology (under review)}, 2021.

\bibitem{me4}
M.~Doff-Sotta, M.~Cannon, and J.~Forbes, ``Spacecraft energy management using
  convex optimisation,'' {\em under review}, 2021.

\bibitem{cvx}
M.~Grant and S.~Boyd, ``{CVX}: Matlab software for disciplined convex
  programming, version 2.1.'' \url{http://cvxr.com/cvx}, Mar. 2014.

\bibitem{mosek}
M.~ApS, {\em Introducing the MOSEK Optimization Suite 9.3.6}, 2021.

\end{thebibliography}
\bibliographystyle{ieeetr}

\end{document}